\documentclass[12pt]{article}
\usepackage{amsmath}
\usepackage{amssymb}
\usepackage{vector}
\font\tencmmib=cmmib10 \skewchar\tencmmib '60
\newfam\cmmibfam
\textfont\cmmibfam=\tencmmib

\def\bbox{\quad\hbox{\vrule \vbox{\hrule \vskip2pt \hbox{\hskip2pt
\vbox{\hsize=1pt}\hskip2pt} \vskip2pt\hrule}\vrule}}
\def\lessim{\ \lower4pt\hbox{$
\buildrel{\displaystyle <}\over\sim$}\ }
\def\gessim{\ \lower4pt\hbox{$\buildrel{\displaystyle >}
\over\sim$}\ }

\def\X{{\cal X}}

\def\la{{\Bigl\langle}}
\def\ra{{\Bigr\rangle}}

\def\qed{\hfill\break\rightline{$\bbox$}}
\newcommand{\e}{\mathbb{E}}

\newcommand{\Reals}{\mathbb{R}}

\newcommand{\Natural}{\mathbb{N}}

\newcommand{\vsi}{{\vec{\sigma}}}

\newtheorem{lemma}{Lemma}
\newtheorem{theorem}{Theorem}
\newtheorem{corollary}{Corollary}
\makeatletter
\@addtoreset{equation}{section}

\makeatother

\font\tencmmib=cmmib10 \skewchar\tencmmib '60
\newfam\cmmibfam
\textfont\cmmibfam=\tencmmib

\def\bbox{\quad\hbox{\vrule \vbox{\hrule \vskip2pt \hbox{\hskip2pt
\vbox{\hsize=1pt}\hskip2pt} \vskip2pt\hrule}\vrule}}
\def\lessim{\ \lower4pt\hbox{$
\buildrel{\displaystyle <}\over\sim$}\ }
\def\gessim{\ \lower4pt\hbox{$\buildrel{\displaystyle >}
\over\sim$}\ }

\def\go0{\to 0}

\def\la{\langle}

\def\leftitem#1{\item{\hbox to\parindent{\enspace#1\hfill}}}

\def\qed{\hfill\break\rightline{$\bbox$}}

\def\ra{\rangle}

\def\sg{\sigma}

\def\sg2{\sigma^2}

\def\__{_{\infty}}
\begin{document}

\title{
Guerra's interpolation using Derrida-Ruelle cascades.}

\author{ 
{\it Dmitry Panchenko}\thanks{Department of Mathematics,
Massachusetts Institute of Technology, Cambridge, MA and
Department of Mathematics, Texas A\&M University, College Station, TX, USA.
Email: panchenk@math.tamu.edu. Partially supported by NSF grant.},
{\it Michel Talagrand}\thanks{
Equipe d'Analyse de l'Institut Math\'ematique, Paris, France
Email: spinglass@talagrand.net.}
}

\date{}

\maketitle

\begin{abstract}
 
New results about Poisson-Dirichlet point processes and Derrida-Ruelle
cascades allow us to  express Guerra's interpolation 
entirely in the language of Derrida-Ruelle cascades and to streamline Guerra's computations.
Moreover, our approach clarifies the nature of the error terms along the interpolation. 
\end{abstract}
\vspace{0.5cm}



\section{Introduction.}

The interpolation invented by Francesco Guerra in \cite{Guerra} 
is one of the most important results in the mathematical theory of 
the Sherrington-Kirkpatrick model \cite{SherK}. Guerra showed for the first time in \cite{Guerra} how the Parisi formula \cite{Parisi} appears naturally as
an upper bound on the free energy. This was a major step toward the 
rigorous proof of this formula in \cite{T-P}. 
One can  define  Guerra's interpolation in terms of Derrida-Ruelle cascades \cite{Ruelle}
similarly to Aizenman-Sims-Starr interpolation \cite{ASS2}; this 
 greatly  simplifies the computation leading to the upper bound on the free energy
(\cite{ASS}, \cite{ASS2}).
However, in order to prove that the upper bound is sharp one needs to
understand precisely the error terms along the interpolation as in \cite{T-P}
(see also \cite{SKgen}) 
and Guerra's original representation is much better suited for this analysis. In this paper we obtain new results about 
Poisson-Dirichlet point processes and Derrida-Ruelle cascades
that allow us to express Guerra's interpolation entirely in the language of the cascades
and, in particular, to easily obtain Guerra's representation of the error terms from
the corresponding representation via Derrida-Ruelle cascades. This interplay not
only streamlines the computations but also helps us understand Guerra's interpolation
on the conceptual level.

We consider a Gaussian Hamiltonian $H_N(\vsi)$ indexed by spin configurations
$\vsi\in\Sigma_N = \{-1,+1\}^N$ with covariance 
\begin{equation}
\e H_N(\vsi^1) H_N(\vsi^2) = \xi(R_{1,2})
\label{cov}
\end{equation}
where 
$$
R_{1,2}=\frac{1}{N}\,\vsi^1\cdot\vsi^2=\frac{1}{N}\sum_{i\leq N}\sigma_i^1 \sigma_i^2
$$ 
is called the overlap of configurations
$\vsi^1,\vsi^2$ and $\xi$ is a smooth convex function such that $\xi(0)=0.$
Given external field parameter $h\in\Reals,$ free energy is defined by 
\begin{equation}
F_N = \frac{1}{N}\e\log \sum_{\vsi} \exp\Bigl(H_N(\vsi)+h\sum_{i\leq N} \sigma_i\Bigr).
\label{FE}
\end{equation}
The external field term $h\sum \sigma_i$ will play
no special role in our considerations so for simplicity of notations
it will be omitted. 

{\it Guerra's interpolation.} Let us first recall Guerra's construction.
Given $k\geq 1,$ consider sequences $\vec{m}$ and  $\vec{q}$
such that
$$
0=m_0<m_1<\ldots<m_{k-1}<m_k=1
$$
and 
$$
0=q_0<q_1<\ldots<q_k<q_{k+1} = 1.
$$
Consider a  matrix 
\begin{equation}
Z=(z_{il}) \,\,\,\mbox{ for }\,\,\, 1\leq i\leq N 
\,\,\,\mbox{ and }\,\,\, 0\leq l\leq k
\label{Z}
\end{equation}
of independent Gaussian r.v. such that 
$
\e z_{il}^2 = \xi'(q_{l+1}) - \xi'(q_l),
$
i.e. the coordinates of each column are i.i.d.
Let 
$$
s=(s_1,\ldots,s_N) 
\,\,\,\mbox{ where }\,\,\, 
s_i=\sum_{0\leq l\leq k} z_{il}.
$$ 
For $0\leq t\leq 1$ we define an interpolating Hamiltonian by
\begin{eqnarray}
&&
H_t(\vsi) 
= 
\sqrt{t} H_N(\vsi) + \sqrt{1-t}\, s\cdot \vsi. 
\label{GuerraH}
\end{eqnarray}
Consider
$
X_{k} = \log \sum_{\vsi} \exp H_t(\vsi)
$
and recursively for $1\leq l\leq k$ define
\begin{equation}
X_{l-1} = \frac{1}{m_l}\log \e_l \exp m_l X_{l}
\label{recX}
\end{equation}
where $\e_l$ denotes the expectation in $(z_{ip})$ for $1\leq i\leq N$ and $l\leq p\leq k.$
By construction, $X_l$ is a function of $(z_{ip})$ for $p\leq l.$ This definition is slightly
different from \cite{T-P}, where $X_l$ denoted what we call $X_{l-1},$ but this indexing
will be more convenient when we define Guerra's interpolation in terms of Derrida-Ruelle
cascades. Finally, we consider
\begin{equation}
\varphi(t) = N^{-1}\,\e X_0.
\label{Guerra}
\end{equation}
It should be obvious that $\varphi(1)=F_N$ and $\varphi(0)$ can be easily computed since all coordinates decouple
and as a result $\varphi(0)$ does not depend on $N.$ 
Let $\theta(x)=x\xi'(x)-\xi(x)$ and for any $a,b\in\Reals$ define
\begin{equation}
\Delta(a,b)=\xi(a)-a\xi'(b) +\theta(b).
\label{Delta}
\end{equation}
By convexity of $\xi,$ $\Delta(a,b)\geq 0.$ The following holds.

\begin{theorem} \label{ThG}{\rm(Guerra)}
We have,
\begin{equation}
\varphi'(t)=
-\frac{1}{2}\,\theta(1)
+\frac{1}{2}\sum_{1\leq r\leq k}(m_r-m_{r-1})\theta(q_r)
-\frac{1}{2}\sum_{1\leq r\leq k} (m_{r} - m_{r-1})
\mu_r\bigl(\Delta(R_{1,2},q_r)\bigr),
\label{phider1}
\end{equation}
where   $\mu_l$ will be described below. 
\end{theorem}

{\it Definition of $\mu_r$}.
Fix $1\leq r\leq k.$
Let 
$$
W_l=\exp m_l (X_{l} - X_{l-1}) \,\,\, \mbox{ for }\,\,\, 1\leq l\leq k.
$$
Notice that by definition of $X_l,$ $W_l$ depends only on $(z_{ip})$ for $p\leq l.$
Consider two copies $Z^1,Z^2$ of $Z$ such that for all $1\leq i\leq N$
\begin{equation}
z_{il}^1=z_{il}^2 
\,\,\,\mbox{ for }\,\,\,
l\leq r-1
\,\,\,\mbox{ and }\,\,\,
z_{il}^1, z_{il}^2
\,\,\,\mbox{ are independent for }\,\,\,
r\leq l.
\label{Z12}
\end{equation}
This means that the columns $0$ through $r-1$ of $Z^1, Z^2$ are completely
correlated and all other columns are independent. We consider Hamiltonians 
$H_t^1$ and $H_t^2$ as above defined in terms of $Z^1$ and $Z^2$ correspondingly
and define $X_l^1, X_l^2$ and $W_l^1, W_l^2$ accordingly. Then, for a function
$f:\Sigma_N^2\to \Reals$ we define
\begin{equation}
\mu_r(f) = 
\e \prod_{1\leq l < r } W_l^1 \prod_{r\leq l\leq k} W_l^1 W_l^2\,
\la f\ra
\label{rem1} 
\end{equation}
where $\la\cdot\ra$ is the Gibbs' average on $\Sigma_N^2$ with respect to
Hamiltonian 
$$
H_t(\vsi^1,\vsi^2)=H_t^1(\vsi^1)+H_t^2(\vsi^2).
$$
Notice that in the first product
for $l<r$ we could also write $W_l^2$ since in this case by construction $W_l^1 = W_l^2.$

{\it Alternative definition of $\mu_r$}.
Fix $1\leq r\leq k.$ Consider a sequence $\vec{n}$ such that 
\begin{equation}
n_l = m_l/2\,\,\,\mbox{ for }\,\,\, l<r
\,\,\,\mbox{ and }\,\,\,
n_l = m_l\,\,\,\mbox{ for }\,\,\, r\leq l.
\label{n}
\end{equation}
In the notations of the first definition let
$
Y_{k} = \log \sum_{\vsi^1,\vsi^2} \exp H_t(\vsi^1,\vsi^2)
$
and recursively for $1\leq l\leq k$ define
$$
Y_{l-1} = \frac{1}{n_l}\log \e_l \exp n_l Y_{l}.
$$
Let $V_l=\exp n_l (Y_{l} - Y_{l-1})$ for $1\leq l\leq k.$
Then, (\ref{rem1}) is equivalent to
\begin{equation}
\mu_r(f) = 
\e \prod_{1\leq l\leq k} V_l\, \la f\ra, 
\end{equation}
where again $\la\cdot\ra$ denotes the Gibbs' average with respect to
the Hamiltonian $H_t(\vsi^1,\vsi^2).$

To see that these definitions are the same, it is a simple exercise to show 
by induction that
$V_l=W_l^1 W_l^2$ for $r\leq l\leq k$ and $V_l=W_l^1 = W_l^2$ for $l<r$
(see Lemma 2.7 in \cite{T-P}).

{\it Guerra's interpolations via Derrida-Ruelle cascades.}
We will now define Guerra's interpolation in the language
of Derrida-Ruelle cascades similarly to \cite{ASS2}.

Given $0<m<1,$ consider a Poisson point process $\Pi$ of intensity measure 
$x^{-1-m}dx$ on $(0,\infty)$. 
Let $(u_n)_{n\geq 1}$ be a decreasing enumeration of $\Pi$
and $w_n=u_n/\sum_{l}u_l.$ The distribution of 
$(w_n)$ is called Poisson-Dirichlet distribution $PD(m,0).$
We will identify a sequence $(u_n)$ with a point process
$\Pi$ and simply call $(u_n)$ itself a Poisson point process.

Let us recall the construction of Derrida-Ruelle cascades 
(see, for example, \cite{Ruelle}, \cite{diluted} or \cite{ASS2})
which involves construction of several processes indexed by $\alpha\in\Natural^k.$
Let us consider a sequence
$$
0<m_1<m_2<\ldots<m_k<1.
$$
We start by constructing a family of point processes on the real line as follows.
\vspace{0.1cm}
\begin{enumerate}
\item[(i)]{
Let $(u_{n_1})_{n_1\geq 1}$ be a decreasing enumeration of 
a Poisson point process on $(0,\infty)$ with intensity measure $x^{-1-m_1}dx.$
}
\item[(ii)]{
Recursively for $2\leq l\leq k,$ for all $(n_1,\ldots, n_{l-1})\in\Natural^{l-1}$
we define independent Poisson point processes
$(u_{n_1\ldots n_{l-1} n_l})_{n_l\geq 1}$
with intensity measure $x^{-1-m_l}dx$ 
independent of all previously constructed processes $(u_{n_1\ldots n_j})$ for $j\leq l-1$.
}
\item[(iii)]{
For $\alpha=(n_1,\ldots,   n_k)\in \mathbb{N}^k$ we define
$
v_{\alpha} = \prod_{1\leq l\leq k} u_{n_1\ldots n_l}
\,\,\,\mbox{ and }\,\,\,
w_{\alpha}=v_{\alpha}\bigr/\sum_{\alpha} v_{\alpha}.
$
}
\end{enumerate}\vspace{0.1cm}
The reason why the sum $\sum v_{\alpha}$ is well defined follows easily from the properties
of Poisson point processes (see, for example, \cite{ASS2}, \cite{diluted}).
We assume that $m_k<1$ is because the sum of Poisson point process corresponding 
to $m_k=1$ is not well defined (equal to $+\infty$ a.s.).
In the interpolation that we will now describe one should formally treat the last step 
corresponding to $m_k=1$ differently but this simple modification will unnecessarily 
complicate the notations. Instead, for simplicity of notations,  
we will work with $m_k<1$ and then formally let $m_k\to 1$.

Let $Z=(z_{0},z_{1},\ldots,z_{k})$ be a column representation of 
a Gaussian matrix in (\ref{Z}). Let us define a sequence $Z_{\alpha}$ of copies of $Z$
as follows.
\begin{enumerate}
\item[(i)]{
Let $(z_{n_1})_{n_1\geq 1}$ be i.i.d. copies of $z_1.$
}
\item[(ii)]{
Recursively for $2\leq l\leq k$, for all $(n_1,\ldots, n_{l-1})\in \Natural^{l-1}$ 
we define independent sequences $(z_{n_1\ldots n_{l-1}n_{l}})_{n_{l}\geq 1}$
of i.i.d. copies of $z_l$ independent of all 
$(z_{n_1\ldots n_{j}})$ for $j\leq l-1.$ 
}
\item[(iii)]{
For all $\alpha=(n_1,\ldots,   n_k)\in\Natural^k$ 
we define $Z_{\alpha}=(z_{il}^\alpha)=(z_0,z_{n_1},z_{n_1 n_2},\ldots,z_{n_1\ldots n_k}).$  
}
\end{enumerate}
Let 
$$
s^{\alpha}=(s_1^{\alpha},\ldots,s_N^{\alpha}) 
\,\,\,\mbox{ where }\,\,\, 
s_i^{\alpha}=\sum_{0\leq l\leq k} z_{il}^{\alpha}.
$$ 
It is easy to check that, by construction, for any $\alpha,\beta\in \Natural^k$
\begin{equation}
\e s_i^{\alpha} s_i^{\beta} = \xi'(q_{\alpha\wedge \beta})
\,\,\,\mbox{ and }\,\,\,
\e s_i^{\alpha} s_j^{\beta} = 0
\,\,\,\mbox{ for }\,\,\, i\not= j
\label{scalpha}
\end{equation}
where 
\begin{equation}
\alpha\wedge\beta = \left\{
\begin{array}{cc}
\min\{l\geq 1: \alpha_l\not = \beta_l\}& \mbox{ if }\,\,\, \alpha\not=\beta\\
k+1 & \mbox{ if }\,\,\, \alpha=\beta.
\end{array}
\right.
\label{wedge}
\end{equation}
For $0\leq t\leq 1$ we define a Hamiltonian
\begin{equation}
H_t(\vsi,\alpha) = 
\sqrt{t} H_N(\vsi) + \sqrt{1-t}\, s^{\alpha}\cdot\vsi
\label{Hta}
\end{equation}
and define
\begin{equation}
\varphi(t)=\frac{1}{N}\e\log \sum_{\alpha,\vsi} w_{\alpha} \exp H_t(\vsi,\alpha).
\label{ASSint}
\end{equation}
Based on the properties of Derrida-Ruelle cascades
we will see that $\varphi(t)$ is equal to Guerra's interpolation
in (\ref{Guerra}). The definition (\ref{ASSint}) is similar to the Aizenman-Sims-Starr interpolation
in \cite{ASS2} with one difference that here we omit an additional term in (\ref{Hta}). In the present setting, 
due to the properties of Derrida-Ruelle cascades, adding this
extra term is a matter of taste. Not adding this term as the advantage to give an interpolation identical
to Guerra's  in (\ref{Guerra}). 
Let us consider a Gibbs probability measure $\Gamma$ on
$\Sigma_N\times \mathbb{N}^{k}$ defined by
\begin{equation}
\Gamma\bigl\{(\vsi,\alpha)\bigr\} \sim w_{\alpha}\exp H_t(\vsi,\alpha).
\label{Gamma}
\end{equation}

\begin{theorem}\label{ThASS}
We have
\begin{equation}
\varphi'(t) =
-\frac{1}{2}\theta(1)
+\frac{1}{2}\e\bigl\la
\theta(q_{\alpha\wedge\beta})
\bigr\ra
-\frac{1}{2}\e\bigl\la
\Delta(R_{1,2},q_{\alpha\wedge\beta})
\bigr\ra
\label{phider2}
\end{equation}
where $\la\cdot\ra$ is the Gibbs average with respect to $\Gamma^{\otimes 2}.$
\end{theorem}
\textbf{Proof.}
By (\ref{ASSint}) and (\ref{Gamma}),
$$
\varphi'(t) 
= 
\frac{1}{2\sqrt{t}} \e\bigl\la H_N(\vsi)\bigr\ra
-
\frac{1}{2\sqrt{1-t}} \e\bigl\la s^{\alpha}\cdot\vsi \bigr\ra.
$$
Using (\ref{cov}) and (\ref{scalpha}), Gaussian integration by parts 
easily implies that this is equal to
\begin{eqnarray*}
\varphi'(t)
&=&
\frac{1}{2}(\xi(1) -\xi'(1))
-\frac{1}{2}
\e\bigl\la \xi(R_{1,2}) - R_{1,2}\xi'(q_{\alpha\wedge\beta})\bigr\ra
\\
&=&
-\frac{1}{2}\theta(1)
+\frac{1}{2}\e\bigl\la
\theta(q_{\alpha\wedge\beta})
\bigr\ra
-\frac{1}{2}\e\bigl\la
\Delta(R_{1,2},q_{\alpha\wedge\beta})
\bigr\ra.
\end{eqnarray*}
and this finishes the proof.
\qed

This proof illustrates that the computation of the derivative in this version
of Guerra's interpolation is a simple exercise compared to the original computation 
of Theorem \ref{ThG} in \cite{Guerra}. However, in Theorem \ref{ThG} the corresponding 
error terms were defined much more precisely and a priori it is not at all obvious
how this can be deduced from (\ref{phider2}). As the following shows, the second term in (\ref{phider2}) is equal to the second
term in (\ref{phider1}).
\begin{theorem}\label{Th6}
For all $1\leq r\leq k$ and for all $0\leq t\leq 1,$
\begin{equation}
\e\la I(\alpha\wedge\beta = r)\ra 
=
\e\Gamma^{\otimes 2}\{\alpha\wedge\beta = r\}
= m_r-m_{r-1}.
\end{equation}
\end{theorem}
This implies that
$$
\e\bigl\la
\theta(q_{\alpha\wedge\beta})
\bigr\ra
=
\sum_{1\leq r\leq k}
\e\bigl\la
I(\alpha\wedge\beta = r)
\bigr\ra\,
\theta(q_{r})
=
\sum_{1\leq r\leq k}(m_r - m_{r-1})\theta(q_r).
$$
It remains to understand the last term in (\ref{phider2}).
Note that in each error term in the last sum in (\ref{phider1}), 
the overlap $R_{1,2}$ is compared to a fixed value $q_{r}.$ 
Therefore, it seems natural that fixing
$\alpha\wedge\beta = r$ in the Gibbs average in (\ref{phider2}) would produce a corresponding
term in (\ref{phider1}). This turns out to be true but the proof will require new
results about Poisson-Dirichlet point processes and Derrida-Ruelle cascades.
\begin{theorem}\label{Th3}
For $1\leq r\leq k,$ we have 
\begin{equation}
\e\bigl\la
\Delta(R_{1,2},q_{\alpha\wedge\beta})I(\alpha\wedge\beta = r)
\bigr\ra
=
(m_r - m_{r-1})\mu_r\bigl(
\Delta(R_{1,2},q_{r})
\bigr).
\label{relate1}
\end{equation}
\end{theorem}
The alternative definition of $\mu_r$ above played an important role in the
proof of Parisi formula in \cite{T-P} and one might be interested in the
corresponding representation via Derrida-Ruelle cascades if one, for example,
wishes to write the interpolation in \cite{T-P} for coupled copies via the cascades.
This can be expressed as follows. Let $(Z^1,Z^2)$ be a pair
of matrices defined in (\ref{Z12}). 
Let $\vec{n}$ be a sequence defined in (\ref{n}) and
let $w_{\alpha}^{(r)}$ be the Derrida-Ruelle cascades corresponding to parameters given by $\vec{n}.$
Next, we generate a sequence $(Z^{1},Z^{2})_{\alpha}$  as above by treating
a pair of matrices as a block matrix with twice as many rows.
We define a Hamiltonian on $\Sigma_N^2\times \mathbb{N}^{k}$ by
\begin{eqnarray}
H_t(\vsi^1,\vsi^2, \alpha) 
&=& 
\sqrt{t} H_N(\vsi^1) + \sqrt{1-t}\, s^{1,\alpha}\cdot\vsi^1 
\nonumber
\\
&+&
\sqrt{t} H_N(\vsi^2) + \sqrt{1-t}\, s^{2,\alpha}\cdot\vsi^2 
\label{Hta2}
\end{eqnarray}
and define a Gibbs' measure $\Gamma_r$ on
$\Sigma_N^2 \times \mathbb{N}^{k}$ by
\begin{equation}
\Gamma_r\bigl\{(\vsi^1,\vsi^2,\alpha)\bigr\} \sim w^{(r)}_{\alpha}\exp H_t(\vsi^1,\vsi^2,\alpha).
\label{Gammar}
\end{equation}
The following holds.
\begin{theorem}\label{Th4} 
For any function $f:\Sigma_N^2\to\Reals$ we have $\mu_r(f) = \e\la f\ra_r$
and, in particular,
\begin{equation}
\mu_r\bigl(
\Delta(R_{1,2},q_{r})
\bigr)
=
\e\la \Delta(R_{1,2},q_{r})\ra_r
\label{mur}
\end{equation}
where $\la\cdot\ra_r$ is the average with respect to
the Gibbs measure $\Gamma_r$ in (\ref{Gammar}).
\end{theorem}

\section{Properties of Poisson-Dirichlet point processes.}\label{SecPD}

%Since
%$$
%\e \card\{x\in \Pi : x\geq 1\} = \int_{1}^{\infty} x^{-1-m}dx =\frac{1}{m} < \infty
%$$
%and
%$$
%\e \sum_{x\in \Pi} xI(x\leq 1) = \int_{0}^{1} x^{-m} dx = \frac{1}{1-m} <\infty,
%$$
%the sum $\sum_{x\in \Pi} x$ is finite a.s. 
In this section we  obtain new results regarding the Poisson-Dirichlet point
process and in the next section we will generalize them to Derrida-Ruelle cascades.
These results will immediately imply Theorems \ref{Th6}, \ref{Th3} and \ref{Th4}.
First, let us state a well-known property of Poisson-Dirichlet point process 
(see \cite{Ruelle} or Lemma 6.5.15 in \cite{SG}).
\begin{lemma}\label{Peq}
Let $0<m<1.$  If $(u_n)$ is a Poisson point process with intensity measure 
$$
d\mu = x^{-1-m}dx
\,\,\,\mbox{ on }\,\, (0,\infty)
$$ 
and $U_n> 0$ are i.i.d. random variables such that $\e U_n <\infty$ then 
$$
(u_n U_n)
\,\,\,\mbox{ and }\,\,\, 
(u_n (\e U_1^m)^{1/m})
$$ 
are both Poisson point processes with the same intensity measure $\e U_1^m d\mu$. 
\end{lemma}
Next, we will prove a result that contains the main idea of the paper.
Let $\X$ be a complete separable metric space that we will also view as
a measurable space with Borel $\sigma$-algebra.
Consider an i.i.d. sequence $(X_n,Y_n)$ with distribution $\nu$ on $\Reals\times \X$ independent 
of $(u_n)$ and such that $X_n>0.$ 
Let $\nu_1,\nu_2$ denote the marginals of $\nu$
and $\nu_x$ denote a regular conditional distribution of $Y$ given $X=x.$
Suppose that $\e X<\infty$ and define by $\nu_{m}$ a probability measure on $\X$
$$
\nu_{m}(B) = \int \frac{x^m}{\e X^m} \nu_x(B)d\nu_1(x)
$$
which is obviously a distribution of $Y$ under the change of density $X^m/\e X^m,$
i.e. for any measurable function $\phi,$
$$
\int \phi(y)d\nu_{m}(y) = \frac{\e X^m \phi(Y)}{\e X^m}.
$$
The following holds.
\begin{lemma}\label{LemNew}
 Poisson point process $(u_n X_n,Y_n)$ has the same distribution as a point
process $((\e X^m)^{1/m}u_n, Y_n')$ where $(Y_n')$ is an i.i.d. sequence independent
of $(u_n)$ with distribution $\nu_{m}.$
\end{lemma}
\textbf{Proof.} 
By the marking theorem (\cite{king}) a point process $(u_n,X_n,Y_n)$ is a Poisson point process
with intensity measure $\mu\otimes \nu$ on $(0,\infty)\times(0,\infty)\times \X.$ By the mapping theorem (\cite{king}), 
$(u_n X_n, Y_n)$ is a Poisson point process with intensity measure given by the image
of $\mu\otimes \nu$ under the mapping $(u,x,y)\to (ux,y)$ if this measure has no atoms.
Let us compute this image measure. Given two measurable sets $A\subseteq (0,\infty)$ and
$B\subseteq \X,$ 
$$
\mu\otimes\nu(ux\in A,y\in B) = 
\int \mu(u: ux\in A)\nu_x(B)d\nu_1(x).
$$
For $x>0$ we have
$$
\mu(u: xu\in A)
= \int I(xu\in A)x^{-1-m}dx = u^m \int I(z\in A)z^{-1-m}dz
= u^m \mu(A)
$$
and, therefore,
$$
\mu\otimes\nu(ux\in A,y\in B) = 
\int x^m \mu(A)\nu_x(B)d\nu_1(x)
=
\e X^m \mu(A) \otimes \nu_{m}(B).
$$
Since measure $\e X^m \mu$ is the intensity measure of a Poisson point process
$((\e X^m)^{1/m} u_n)$ this finishes the proof.
\qed

As an application of Lemma \ref{LemNew} we will give a new simple proof
of Theorem 6.4.5 in \cite{SG}.
\begin{corollary}\label{Cor1}
If $(X_n,Y_n)$ are i.i.d. such that $X\geq 1$ and $\e X^2, \e Y^2 <\infty$ then
\begin{eqnarray}
&&
\hspace{-1.5cm}\e \frac{\sum u_n Y_n}{ \sum u_n X_n} 
=
\frac{\e X^{m-1}Y}{ \e X^m} ,
\\
&&
\e \frac{\sum u_n^2 Y_n^2}{ (\sum u_n X_n)^2} 
=
(1-m)\frac{\e X^{m-2} Y^2}{ \e X^m} ,
\label{t2}
\\
&&
\hspace{1.5cm}\e \frac{\sum_{n\not=m} u_n u_m Y_n Y_m}{ (\sum u_n X_n)^2} 
=
m\Bigl(\frac{\e X^{m-1} Y}{ \e X^m}\Bigr)^2.
\end{eqnarray}
\end{corollary}
\textbf{Proof.}
If we denote by $c=(\e X^m)^{1/m}$ then by Lemma \ref{LemNew},
$$
\e \frac{\sum u_n Y_n}{ \sum u_n X_n} 
=
\e \frac{\sum (u_n X_n) (Y_n/X_n)}{ \sum u_n X_n} 
= 
\e \frac{\sum (u_n c) (Y_n/X_n)'}{ \sum u_n c} 
=
\e \frac{X^m}{ \e X^m} \frac{Y}{X}
$$
since the markings $(Y_n/X_n)'$ are independent of $(u_n)$
and the distribution is given by the change of density $X^m/\e X^m.$
Similarly,
\begin{eqnarray*}
\e \frac{\sum u_n^2 Y_n^2}{ (\sum u_n X_n)^2} 
&=&
\e \frac{\sum (u_n X_n)^2 (Y_n/X_n)^2}{( \sum u_n X_n)^2} 
\\
&=&
\e \frac{\sum (u_n c)^2 (Y_n/X_n)'^2}{ (\sum u_n c)^2} 
=
\e \frac{X^m}{ \e X^m} \frac{Y^2}{X^2}\,
\e \frac{\sum u_n^2 }{ (\sum u_n)^2}. 
\end{eqnarray*}
To finish the proof of (\ref{t2}) it remains to use a well-known fact 
(Corollary 2.2 in \cite{Ruelle} or Proposition 1.2.7 in \cite{SG})
\begin{equation}
\e \sum w_n^2 =(1-m).
\label{rue}
\end{equation}
Finally,
\begin{eqnarray*}
\e \frac{\sum_{n\not=m} u_n u_m Y_n Y_m}{ (\sum u_n X_n)^2} 
&=&
\e \frac{\sum_{n\not= m} (u_n X_n) (u_m X_m) (Y_n/X_n)(Y_m/X_m)}{( \sum u_n X_n)^2} 
\\
&=& 
\e \frac{\sum_{n\not = m} (u_n c) (u_m c) (Y_n/X_n)'(Y_m/X_m)'}{ (\sum u_n c)^2} 
=
m\Bigl(\e\frac{X^m}{ \e X^m} \frac{Y}{X}\Bigr)^2
\end{eqnarray*}
since by (\ref{rue}),
$\e \sum_{n\not=m} w_n w_m  = 1-\e \sum w_n^2 =m.$
\qed

\section{Properties of Derrida-Ruelle cascades.}\label{SecDR}

%Let us start with a comment about a Derrida-Ruelle process 
%$(u_{n_1},\ldots,u_{n_1\ldots n_k})$ index by $\alpha=(n_1,\ldots,   n_k)\in\Natural^k$
%defined in the beginning of Section \ref{ASSI}. As in the previous section,
%we use index $\alpha$ for convenience of notations only and assume that any other
%process indexed by $\alpha$ is defined conditionally on the Derrida-Ruelle
%process and is indexed by it. 

Let us construct a general random process $Z_{\alpha}$ indexed by $\alpha\in \mathbb{N}^k$  
in a much more general way than the random matrix process in the second version of
Guerra's interpolation above.
Consider complete separable metric spaces $\X_1,\ldots,\X_k$ which we also view
as measurable spaces with Borel $\sigma$-algebras and for $1\leq l\leq k$ let
$$
\X^l = \X_1\times\ldots\times \X_l.
$$
Consider a probability measure $\nu$ on $\X_1$ and for $1\leq l< k$
consider regular conditional distributions 
\begin{equation}
\nu_l(\cdot|x) \,\,\,\mbox{ on }\,\,\, \X_{l+1} \,\,\,\mbox{ for }\,\, x\in\X^l.
\label{condd}
\end{equation}
We generate a process
$$
Z_{\alpha}=(z_{n_1},z_{n_1 n_2},\ldots,z_{n_1 n_2\ldots n_k})\in{\cal X}^k
$$
according to the following recursive procedure.
\begin{enumerate}
\item[(i)]{
Generate i.i.d. random variables $(z_{n_1})_{n_1\geq 1}$  with distribution $\nu.$ 
}
\item[(ii)]{
Recursively over $2\leq l\leq k$, given $(z_{n_1},\ldots,z_{n_1\ldots n_{l-1}})$
for all $n_1\ldots n_{l-1}\in \Natural$ 
we generate i.i.d. sequences $(z_{n_1\ldots n_{l-1} n_{l}})_{n_{l}\geq 1}$ 
with distributions 
\begin{equation}
\nu_l(\cdot|z_{n_1},\ldots,z_{n_1\ldots n_{l-1}})
\label{recD1}
\end{equation}
independently for all $n_{1},\ldots,n_{l-1}.$
}
\item[(iii)]{
For each $\alpha=(n_1,\ldots,   n_k)\in\Natural^k$ 
we define $Z_{\alpha}=(z_{n_1},z_{n_1 n_2},\ldots,z_{n_1\ldots n_k})\in{\cal X}^k.$  
}
\end{enumerate}
For convenience of notations, given $\alpha=(n_1,\ldots,   n_k)$ we denote
for $1\leq l\leq k,$
\begin{equation}
\alpha^l = (n_1\ldots n_l),\,\,\,\,\,
u_{\alpha^l}=u_{n_1\ldots n_l}
\,\,\,\mbox{ and }\,\,\,
v_{\alpha^l} = \prod_{1\leq j\leq l} u_{\alpha^j}
\label{uv}
\end{equation}
so that $v_{\alpha^{l+1}}=v_{\alpha^l} u_{\alpha^l}.$
Given $Z_{\alpha}\in{\cal X}^k$  we denote 
$$
z_{\alpha^l} = z_{n_1\ldots n_l}
\,\,\,\mbox{ and }\,\,\,
Z_{\alpha^l} = 
(z_{n_1},\ldots,z_{n_1\ldots n_l}).
$$
Consider a measurable function $X:\X^k\to\Reals$ such that 
$\e \exp X(Z_{\alpha}) <\infty.$
Let $X_{\alpha}=X(Z_{\alpha})$ and recursively for $1\leq l\leq k$ define
\begin{equation}
X_{\alpha^{l-1}}=\frac{1}{m_{l}}\log \e_{l} \exp m_{l} X_{\alpha^{l}}
\label{Xrec}
\end{equation}
where $\e_l$ denotes the expectation conditionally on 
$(Z_{\alpha^{l-1}})_{\alpha\in\Natural^k}$ and
\begin{equation}
W_{\alpha^l}=\exp m_l(X_{\alpha^l} - X_{\alpha^{l-1}}). 
\label{Wl}
\end{equation}
Thus, both $X_{\alpha^l}$ and $W_{\alpha^l}$ are functions of
$Z_{\alpha^l}.$ 
In particular, $X_0:= X_{\alpha^0}$ is a constant.
It will be convenient to think of $W_{\alpha^l}$
as a function of two variables
$$
W_{\alpha^l} = W_{l}(Z_{\alpha^{l-1}},z_{\alpha^l}).
$$
Let us now generate another process $Z_{\alpha}'$ exactly the same way
as $Z_{\alpha}$ with one modification that instead of (\ref{recD1})
the distribution of $(z_{n_1\ldots n_{l-1} n_{l}}')_{n_{l}\geq 1}$ conditionally on
$Z_{\alpha^{l-1}}' = (z_{n_1}',\ldots,z_{n_1\ldots n_{l-1}}')$ will be given
by
\begin{equation}
W_{l}(Z_{\alpha^{l-1}}', x)\, d\nu_l(x|Z_{\alpha^{l-1}}').
\label{DC}
\end{equation}
This is a probability measure because by (\ref{Xrec}), (\ref{Wl}) and (\ref{recD1}),
$$
\int W_{l}(Z_{\alpha^{l-1}}', x)\,d\nu_l(x|Z_{\alpha^{l-1}}')
=
\e_l \exp m_l(X_{\alpha^l} - X_{\alpha^{l-1}})
=1.
$$
For $1\leq l\leq k,$ let us define
\begin{equation}
e_{\alpha^l}=
\exp (X_{\alpha^l} - X_{\alpha^{l-1}}).
\label{El}
\end{equation}
The following in the generalization of Lemma \ref{LemNew}.

\begin{lemma}\label{Lem3}
The point processes 
\begin{equation}
(u_{\alpha^1}e_{\alpha^1},\ldots,u_{\alpha^k} e_{\alpha^k},Z_{\alpha^k})
\,\,\,\mbox{ and }\,\,\, 
(u_{\alpha^1},\ldots,u_{\alpha^k}, Z_{\alpha^k}')
\label{Lem41}
\end{equation}
on ${\Reals^+}^{k}\times {\cal X}^k$ have the same distribution.
\end{lemma}
\textbf{Proof.}
The proof is by induction on $k.$ The case $k=1$ immediately follows from Lemma \ref{LemNew}. 
Consider $k>1.$ By induction assumption, point processes
\begin{equation}
(u_{\alpha^1}e_{\alpha^1},\ldots,u_{\alpha^{k-1}} e_{\alpha^{k-1}},Z_{\alpha^{k-1}})
\,\,\,\mbox{ and }\,\,\, 
(u_{\alpha^1},\ldots,u_{\alpha^{k-1}}, Z_{\alpha^{k-1}}')
\label{indu}
\end{equation}
have the same distribution. 
If we write 
$$
Z_{\alpha^k}=(Z_{\alpha^{k-1}},z_{\alpha^k})
\,\,\,\mbox{ and }\,\,\,
Z_{\alpha^k}'=(Z_{\alpha^{k-1}}',z_{\alpha^k}')
$$
it suffices to show that conditionally on the processes (\ref{indu}), the two processes
\begin{equation}
\bigl(u_{\alpha^k} e_{\alpha^{k}},z_{\alpha^k}\bigr)
\,\,\,\mbox{ and }\,\,\, 
\bigl(u_{\alpha^k},z_{\alpha^k}'\bigr)
\label{ppk1}
\end{equation}
have the same distribution .
Let us write $\alpha^k=(\alpha^{k-1},n)$ and for a fixed $\alpha^{k-1}$
look at the point process
$\bigl(u_{\alpha^k} e_{\alpha^{k}},z_{\alpha^k}\bigr)_{n\geq 1}.$
Let us apply Lemma \ref{LemNew} to this sequence conditionally on (\ref{indu}). 
By (\ref{Xrec}), 
$$
\e_k e_{\alpha^k}^{m_k}
=
\e_k \exp m_k (X_{(\alpha^{k-1},n)} - X_{\alpha^{k-1}}) = 1
$$
and, therefore, by Lemma \ref{LemNew}, the  point processes
\begin{equation}
\bigl(u_{\alpha^k} e_{\alpha^{k}},z_{\alpha^k}\bigr)_{n\geq 1}
\,\,\,\mbox{ and }\,\,\, 
\bigl(
u_{(\alpha^{k-1},n)},
z_{(\alpha^{k-1},n)}'
\bigr)_{n\geq 1}
\label{ppk2}
\end{equation}
have the same distribution, 
where $z_{(\alpha^{k-1},n)}'$ is distributed as $z_{(\alpha^{k-1},n)}$
under the change of density 
$$
\frac{e_{\alpha^k}^{m_k}}
{\e_k e_{\alpha^k}^{m_k}}
=
\exp m_k (X_{(\alpha^{k-1},n)} - X_{\alpha^{k-1}})
=
W_k(Z_{\alpha^{k-1}},z_{\alpha^k}).
$$
By construction, $z_{(\alpha^{k-1},n)}$ are distributed according to  
$\nu_l(\cdot|Z_{\alpha^{k-1}})$ and the change of density defines a distribution
$$
W_{k}(Z_{\alpha^{k-1}}, x)\,d\nu_k(x|Z_{\alpha^{k-1}})
$$
which is precisely the distribution (\ref{DC}) for $l=k$.  
Since conditionally on (\ref{indu}) processes (\ref{ppk2}) are generated
independently for all $\alpha^{k-1},$ this shows that conditionally on (\ref{indu})
both processes in (\ref{ppk1}) are generated according to the same distribution 
and this finishes the proof.
\qed

In particular, Lemma \ref{Lem3} implies that the processes
\begin{equation}
v_{\alpha} \exp (X_{\alpha}-X_0)
=
\prod_{1\leq l\leq k} u_{\alpha^l} e_{\alpha^l} 
\,\,\,\mbox{ and }\,\,\,
v_{\alpha}=\prod_{1\leq l\leq k} u_{\alpha^l}
\label{Lem3c}
\end{equation}
have the same distribution, which generalizes Theorem 5.4 in \cite{ASS2}.
As a consequence we get (Proposition 2 in \cite{diluted}) 
\begin{equation}
\e\log \sum w_{\alpha}\exp X_\alpha = X_0.
\label{DRlog}
\end{equation}
Using (\ref{DRlog}) one only needs to compare the definitions to observe 
the equality of (\ref{Guerra}) and (\ref{ASSint}). 
Using (\ref{Lem3c}), Lemma \ref{Lem3} also implies that
\begin{equation}
\bigl(v_{\alpha} \exp (X_{\alpha}-X_0),Z_{\alpha}\bigr)
\,\,\,\mbox{ and }\,\,\,
\bigl(v_{\alpha}, Z_{\alpha}'\bigr)
\label{Lem3cor}
\end{equation}
have the same distribution. 
As we will now show, this immediately implies Theorems \ref{Th3} and \ref{Th4}.
Moreover, the change of density (\ref{DC}) makes the definition of measures 
$\mu_r$ in Guerra's interpolation in (\ref{phider1}) much more transparent. 

In addition to $X$, consider a measurable function
$Y:\X^k\to\Reals$ such that $\e Y^2(Z_{\alpha}) <\infty$
and let $Y_{\alpha} = Y(Z_{\alpha}).$ 
Theorem \ref{Th4} is an immediate consequence of the following.
\begin{theorem}\label{Th7}
We have
\begin{equation}
\e\frac{\sum_{\alpha} v_{\alpha} (\exp X_{\alpha})\, Y_{\alpha}}
{\sum_{\alpha} v_{\alpha} \exp X_{\alpha}}
=
\e \prod_{1\leq l\leq k} W_{\alpha^l} Y_{\alpha}.
\label{Gerror2}
\end{equation}
\end{theorem}
\textbf{Proof.}
The proof follows immediately by (\ref{Lem3cor}), because
\begin{eqnarray*}
\e\frac{\sum_{\alpha} v_{\alpha} (\exp X_{\alpha})\, Y(Z_{\alpha})}
{\sum_{\alpha} v_{\alpha} \exp X_{\alpha}}
&=&
\e \sum w_{\alpha} Y(Z_{\alpha}')
\\
&=&
\e Y(Z_{\alpha}') 
=
\e \prod_{1\leq l\leq k} W_{\alpha^l} Y_{\alpha},
\end{eqnarray*}
where in the second line $\alpha$ is fixed and the last
equality holds since the distribution of $Z_{\alpha}'$ 
is defined by the change of density (\ref{DC}).
\qed

Let us now fix $1\leq r\leq k.$ Consider a measurable function  
$
Y:\X^k\times \X^k \to \Reals
$
such that $\e Y^2(Z_{\alpha},Z_{\beta})<\infty$ for any $\alpha,\beta\in\Natural^k$ 
and let $Y_{\alpha,\beta} = Y(Z_{\alpha},Z_{\beta}).$
Let us consider fixed $\alpha,\beta\in \mathbb{N}^k$
such that $\alpha\wedge\beta=r.$
Let
$$
M_r 
=
\e \prod_{l< r} W_{\alpha^l} \prod_{l\geq r} W_{\alpha^{l}} W_{\beta^{l}} Y_{\alpha,\beta}.
$$
Clearly, $M_r$ depends on $\alpha$ and $\beta$ only through $r=\alpha\wedge\beta.$
Theorem \ref{Th3} is an immediate consequence of the following.
\begin{theorem}
We have
\begin{equation}
\e\frac{\sum_{\alpha\wedge\beta=r} v_{\alpha} v_{\beta} \exp(X_{\alpha}+X_{\beta})\, 
Y_{\alpha,\beta}}
{\left(\sum_{\alpha} v_{\alpha} \exp X_{\alpha}\right)^2}
=
(m_{r}-m_{r-1})M_r.
\label{Gerror}
\end{equation}
\end{theorem}
\textbf{Proof.}
Again, by (\ref{Lem3cor})
$$
\e\frac{\sum_{\alpha\wedge\beta=r} v_{\alpha} v_{\beta} \exp(X_{\alpha}+X_{\beta})\, 
Y_{\alpha,\beta}}
{\left(\sum_{\alpha} v_{\alpha} \exp X_{\alpha}\right)^2}
=
\e Y(Z_{\alpha}',Z_{\beta}')\,
\e \sum_{\alpha\wedge\beta = r} w_{\alpha}w_{\beta},
$$
where $\e Y(Z_{\alpha}',Z_{\beta}')$ is taken for any fixed $\alpha$ and $\beta$
such that $\alpha\wedge\beta=r.$ 
By construction, this expectation is equal to $M_r$ because the distribution
of $Z_{\alpha}'$ is defined by the change of density (\ref{DC})
and, because, since $\alpha\wedge\beta = r,$ the function $Y(Z_{\alpha}',Z_{\beta}')$
depends on one copy $z_{\alpha^l}'=z_{\beta^l}'$  for $l<r$ and on two independent
copies $z_{\alpha^l}'$ and $z_{\beta^l}'$ for $l\geq r.$
It remains to show that 
\begin{equation}
\e \sum_{\alpha\wedge\beta = r} w_{\alpha}w_{\beta}
= m_r - m_{r-1}.
\label{thproo}  
\end{equation}
Given $\alpha\in\Natural^k$ let us write $\alpha^r=(a,n)$
for $a\in\Natural^{r-1}$ and $n\in \Natural.$ If $\alpha\wedge\beta=r$ then $\beta^r=(a,m)$
for $m\not = n.$ In the notations of (\ref{uv}) let us define
$U_{(a,n)} = \sum_{\gamma : \gamma^r = (a,n)} \prod_{r<l\leq k} u_{\gamma^l}.$
Then
$$
\sum_{\alpha\wedge\beta = r} w_{\alpha}w_{\beta}
=
\frac{\sum_{\alpha\wedge\beta=r} v_{\alpha}v_{\beta}}  
{\bigl(\sum_{\alpha} v_{\alpha}\bigr)^2}
=
\frac{\sum_{a} v_a^2 \sum_{n\not=m} (u_{(a,n)} U_{(a,n)}) (u_{(a,m)} U_{(a,m)})}
{\bigl(\sum_{a} v_a \sum_{n} u_{(a,n)} U_{(a,n)}\bigr)^2}.
$$
A sequence $(U_{(a,n)})$ is i.i.d. by construction and, therefore, by Lemma \ref{Peq},
a point process $(u_{(a,n)} U_{(a,n)})$ has the same distribution as
$(u_{(a,n)} c)$ where $c=(\e U_{(a,n)}^{m_r})^{1/m_r}<\infty.$
As a result,
$$
\e \sum_{\alpha\wedge\beta = r} w_{\alpha}w_{\beta}
=
\e \frac{\sum_{a} v_a^2 \sum_{n\not=m} u_{(a,n)} u_{(a,m)}}
{\bigl(\sum_{a} v_a \sum_{n} u_{(a,n)}\bigr)^2}.
$$
Using that
$$
\sum_{n\not= m}u_{(a,n)}u_{(a,m)} = \bigl(\sum_{n}u_{(a,n)}\bigr)^2 - \sum_{n} u_{(a,n)}^2
=
U_a^2 - \sum_{n} u_{(a,n)}^2
$$ 
where we introduced $U_a = \sum_{n}u_{(a,n)},$ we can write
\begin{equation}
\e\frac{\sum_{a} v_{a}^2 \sum_{n\not =m} u_{(a,n)} u_{(a,m)} }
{\left(\sum_{a} v_{a} \sum_{n} u_{(a,n)}\right)^2}
=
\e\frac{\sum_{a} (v_{a} U_{a})^2 }
{\left(\sum_{a} v_{a} U_a\right)^2}
-
\e\frac{\sum_{a,n} v_{(a,n)}^2  }
{\bigl(\sum_{a,n} v_{(a,n)}\bigr)^2}.
\label{diffm}
\end{equation}
By Corollary 3.3 in \cite{Ruelle}, the process $(v_{a,n}/\sum v_{a,n})$ has
Poisson-Dirichlet distribution $PD(m_r,0).$ By Lemma \ref{Lem3} above,
the process $(v_a U_a/\sum v_a U_a)$ has the same distribution as the process
$(v_a/\sum v_a)$ which again, by Corollary 3.3 in \cite{Ruelle}, is $PD(m_{r-1},0).$
Therefore, using (\ref{rue}) twice implies that the right hand side of (\ref{diffm})
is equal to $(1-m_r)-(1-m_{r-1}) = m_r-m_{r-1}.$ 
This finishes the proof.
\qed

Finally, we prove Theorem \ref{Th6}.

\textbf{Proof of Theorem \ref{Th6}.}
Let $\Gamma_1$ be a marginal on $\Natural^k$
of measure $\Gamma$ defined in (\ref{Gamma}). Then 
\begin{equation}
\Gamma_1\{\alpha\} = v_{\alpha} f_{\alpha}\Bigr/ \sum_{\alpha} v_{\alpha} f_{\alpha}
\,\,\,\mbox{ where }\,\,\,
f_{\alpha} = \sum_{\vsi} \exp H_t(\vsi,\alpha).
\label{Gamma1}
\end{equation}
By Lemma \ref{Lem3}, conditionally on $H_N(\vsi)$ and $(z_{i0}, y_{i0})_{1\leq i\leq N},$
the sequence $(\Gamma_1\{\alpha\})_{\alpha\in\Natural^k}$ is equal in distribution to
the sequence $(w_{\alpha})_{\alpha\in\Natural^k}$ and, consequently, the same is true unconditionally.
Therefore,
$$
\e\Gamma^{\otimes 2}\{\alpha\wedge\beta = r\}
=
\e\Gamma_1^{\otimes 2}\{\alpha\wedge\beta = r\}
=
\e \sum_{\alpha\wedge\beta = r} w_{\alpha} w_{\beta}
=
m_r - m_{r-1},
$$
using (\ref{thproo}). This finishes the proof.
\qed

\end{document}